\documentstyle[12pt]{article}
\newtheorem{Lemma}{Lemma}

\newtheorem{Prop}{Proposition}

\newtheorem{Cor}{Corollary}

\newtheorem{Conj}{Conjecture}
\newcommand{\pf}{\medskip\noindent{\sc Proof: }}
\newcommand{\qed}{$\Box$}
\newcommand{\rTo}{\,\longrightarrow\,}

\newcommand{\DS}{\displaystyle}
\begin{document}
\title{Factorials and powers, a minimality result, revisited}
\author{David E. Radford \\
Department of Mathematics, Statistics \\
and Computer Science (m/c 249)    \\
851 South Morgan Street   \\
University of Illinois at Chicago\\
Chicago, Illinois 60607-7045
 }
\maketitle
\date{}

\begin{abstract}
{\small  \rm Let $a > 1$. Then $a^n < n!$ for some positive integer $n$. There are several numerical sequences associated with the study of the smallest such integer which are studied in \cite{RadFact} and \cite{RadGamma}. Here we continue the examination of one of them.}
\end{abstract}
\let\thefootnote\relax\nonumber
%
%SECTION
%
\section{Introduction}
Let $a > 1$ be a fixed real number. Then $a^\ell \leq \ell!$ for some positive integer $\ell$. Let $n_a$ be the smallest such integer. Now $\DS{\frac{n}{e} < a \leq \frac{n+1}{e}}$ for some $n \geq 2$. Let $\sigma_n$ be the largest integer $\ell$ such that $n + (\ell - 1) \leq e\sqrt[n]{n!}\;$. Then $n_a = n - \sigma_n + r$, where $r \in \{ 1, 2, 3 \}$, by \cite[Theorem 3]{RadFact}.

The sequence $\sigma_1, \sigma_2, \sigma_3 \ldots \;$ is increasing and seems to increase very slowly. For example, $\sigma_n = 15$ when $n = 10^{12}$. This and other examples are found in the discussion at the end of \cite[\S 6]{RadFact}. In this paper we show how to narrow $\sigma_n$ to one of two possibilities when $n \geq 4$ using very simple estimates. We also examine the growth rate of the sequence.
\section{A basic sequence associated with factorials}\label{SecSigmaN}
In this section we draw heavily from the material of \cite{RadFact}. First we describe basic functions and some of the results found therein which are needed for this paper.

Let $x > 0$. These basic functions which we use here are the four
$$
L(x) = \displaystyle{\pi^\frac{1}{2x}}e^{\frac{1}{(12x+1)x}}, \;\; R(x) = \pi^\frac{1}{2x}e^\frac{1}{12x^2}, \;\; P(x) = x^{\frac{1}{x}}, \;\; \mbox{and} \;\; T(x) = xP(x)
$$
which are positive valued and infinitely differentiable. For a function $F(x)$ set $a_F(x) = (F(x) - 1)x$. Then $\DS{F(x) = 1 + \frac{a_F(x)}{x}}$. It is easy to see that $\lim_{x \rTo \infty}F(x) = 1$, where $F(x)$ is any one of the first three of the four functions mentioned at the beginning of this paragraph. For any one of these three, limit statement suggests that it is natural to study $F(x)$ in terms of $a_F(x)$.

Let $n \geq 1$. An important estimate for us is Robbins' improvement \cite{Robbins} on Sterling's approximation of factorials. It is equivalent to
$$
\frac{1}{2}L(n)T(2n) < e\sqrt[n]{n!} < \frac{1}{2}R(n)T(2n)
$$
which is (2) of \cite{RadFact}. The previous expression is a very minor reformulation of Robbins' result. The form of Robbins' result most useful to us here is
\begin{equation}\label{EqA}
L(n)P(2n) < e\left(\frac{\sqrt[n]{n!}}{n}\right) < R(n)P(2n).
\end{equation}

The results of \cite[\S 5]{RadFact} apply to the sequence $T_1, T_2, T_3, \ldots\;$ whose terms are $T_i = e\sqrt[i]{i!}\;$ for all $i \geq 1$. Therefore, as $\sigma_1 = 2$,
\begin{equation}\label{EqB}
2 = \sigma_1 \leq \sigma_2 \leq \sigma_3 \leq \cdots
\end{equation}
by \cite[(26)]{RadFact} and
\begin{equation}\label{EqC}
\sigma_{i+1} = \sigma_i \;\;\ \mbox{or} \;\;  \sigma_{i+1} = \sigma_i + 1
\end{equation}
for all $i \geq 1$ by \cite[(27)]{RadFact}.

As $n + \sigma_n - 1 \leq e\sqrt[n]{n!} < n+ \sigma_n$, it follows that
$$
\frac{\sigma_n}{n} - \frac{1}{n} \leq e\left(\frac{\sqrt[n]{n!}}{n}\right) - 1 < \frac{\sigma_n}{n}.
$$
These inequalities yield
\begin{equation}\label{EqD}
L(n)P(2n) - 1 < \frac{\sigma_n}{n} < R(n)P(2n) - 1 + \frac{1}{n}
\end{equation}
by (\ref{EqA}). Set ${\mathcal L}(x) = L(x)P(2x)$ and ${\mathcal R}(x) = L(x)P(2x)$. Then (\ref{EqD}) can be reformulated
\begin{equation}\label{EqDD}
a_{{\mathcal L}(n)} < \sigma_n < a_{{\mathcal R}(n)} + 1.
\end{equation}
For functions $f(x)$ and $g(x)$, note that
$$
\DS{a_{fg}(x) = \frac{a_f(x)a_g(x)}{x} + a_f(x) + a_g(x)} = a_f(x)g(x) + a_g(x).
$$
Let $p(x) = P(2x)$. Then
\begin{equation}\label{EqEEE}
a_{\mathcal L}(x) - a_p(x) = a_L(x)P(2x),
\end{equation}
\begin{equation}\label{EqFFF}
a_{\mathcal R}(x) - a_p(x) = a_R(x)P(2x),
\end{equation}
and therefore
\begin{equation}\label{EqDDD}
a_{\mathcal R}(x) - a_{\mathcal L}(x) = (a_R(x)- a_L(x))P(2x).
\end{equation}
\begin{Lemma}\label{LemmaSigmaXEstimates}
Let $x > 0$. Then:
\begin{enumerate}
\item[{\rm (a)}] $a_{\mathcal R}(x) > a_{\mathcal L}(x) > a_p(x) > 0$.
\item[{\rm (b)}] $a_{\mathcal L}(x)- a_p(x)$ and $a_{\mathcal R}(x) - a_p(x)$ are strictly decreasing functions on $[e/2, \infty)$.
\item[{\rm (c)}] $\DS{\lim_{x \rTo \infty}(a_{\mathcal L}(x)- a_p(x)) = \lim_{x \rTo \infty}(a_{\mathcal R}(x)- a_p(x)) = \frac{\ln\pi}{2}}$.
\item[{\rm (d)}] $\lim_{x \rTo \infty} (a_{\mathcal R}(x) - a_{\mathcal L}(x)) = 0$.
\end{enumerate}
\end{Lemma}

\pf
First of all $R(x) > L(x) > 1$. This means $a_R(x) > a_L(x) > 0$. Since $p(x) > 1$ also, part (a) follows by (\ref{EqEEE}) and (\ref{EqDDD}). Now $a_L(x)$ and $a_R(x)$ are strictly decreasing functions on $((\ln\pi/2), \infty)$ by \cite[Proposition 2(a)]{RadFact}. It was shown in \cite[\S 1]{RadFact} that $P(x)$ is strictly decreasing on $[e, \infty)$. Part (b) now follows by (\ref{EqEEE}) and (\ref{EqFFF}). Since $\lim_{x \rTo \infty} a_L(x) = \lim_{x \rTo \infty}a_R(x) = \DS{\frac{\ln \pi}{2}}$, which follows by \cite[Proposition 2(b)]{RadFact}, and $\lim_{x \rTo \infty}P(x) = 1$, part (c) now follows by (\ref{EqEEE}) and (\ref{EqFFF}) and part (d) by (\ref{EqDDD}).
\qed
\medskip

Let $n_1 < n_2 < n_3 < \ldots \;$ be the sequence defined by $\sigma_1 = \cdots = \sigma_1$ and $\sigma_{n_i+1} = \sigma_{n_i} + 1 = \cdots = \sigma_{n_{i+1}}$ for $i \geq 1$. Then $n_1, n_2, \ldots \;$ mark the positions in $\sigma_1, \sigma_2, \sigma_3, \ldots \;$ just before changes of term value. We are interested in the ratio $\DS{\frac{n_{i+1}}{n_i}}$.

Set $\delta(x) = a_{\mathcal R}(x) - a_{\mathcal L}(x)$, let $i \geq 1$, and set $d(i) = \max\{\delta(n_{i+1}), \delta(n_i)\}$. Then $d(i) > 0$ by part (a) of Lemma \ref{LemmaSigmaXEstimates} and $\lim_{i \rTo \infty}d(i) = 0$ by part (d) of the same. Let $c = \sigma_{n_{i+1}} = \sigma_{n_i} + 1$. Then
$$
a_{\mathcal L}(n_i) + 1 < c < a_{\mathcal R}(n_i) + 2 \;\; \mbox{and} \;\; a_{\mathcal L}(n_{i+1}) < c < a_{\mathcal R}(n_{i+1}) + 1
$$
by (\ref{EqDD}). Therefore $|a_{\mathcal L}(n_{i+1}) - (a_{\mathcal L}(n_{i})+1)| < d(i)$, or equivalently
\begin{equation}\label{EqFFFF}
1 - d(i) < a_{\mathcal L}(n_{i+1}) - a_{\mathcal L}(n_{i}) < 1 + d(i).
\end{equation}
Thus the sequence of positive terms $a_{\mathcal L}(n_1) < a_{\mathcal L}(n_2) < a_{\mathcal L}(n_3) < \ldots \;$ is eventually strictly increasing,
\begin{equation}\label{EqF1FF}
\lim_{i \rTo \infty} (a_{\mathcal L}(n_{i+1}) - a_{\mathcal L}(n_{i})) = 1,
\end{equation}
and thus
\begin{equation}\label{EqFXFF}
\lim_{i \rTo \infty}\left(\frac{a_{\mathcal L}(n_{i+1})}{a_{\mathcal L}(n_{i})}\right) = 1.
\end{equation}

We find estimates for $a_{\mathcal L}(x)$ and $a_{\mathcal R}(x)$ which connect them to $\ln (2x)$. At this point we assume $\DS{x \geq \frac{e}{2}}$. Set $\DS{\alpha = \frac{\ln \pi}{2}}$. Observe that $\DS{\frac{e}{2} > \alpha > \frac{1}{2}}$. Now $\DS{R(x) < 1 + \frac{1}{x}}$ for all $\DS{x \geq \frac{e}{2}}$ by  \cite[Corollary 1]{RadFact} as $x > 1$. Next we note $\DS{P(2x) < \frac{2x}{2x - \ln (2x)} = 1 + \frac{\ln (2x)}{2x - \ln(2x)}}$ which follows by \cite[Lemma 1]{RadFact} as $2x > 1$. Therefore
$$
R(x)P(2x) - 1 < \left(1 + \frac{1}{x}\right)\left(1 + \frac{\ln (2x)}{2x - \ln(2x)}\right) - 1 = \frac{2x + x\ln(2x)}{x(2x - \ln(2x))}
$$
which means that $a_{\mathcal R}(x) < \ln(2x)QR(x)$, where
$$
QR(x) = \left(\frac{2}{\ln(2x)} + 1\right)\left(\frac{x}{2x - \ln(2x)}\right).
$$
Now $QR(x)$ is the product of two positive valued strictly decreasing functions, since $\DS{x > \frac{e}{2}}$; therefore $QR(x)$ is of the same type. It is easy to see that $\lim_{x \rTo \infty}QR(x) = \DS{\frac{1}{2}}$.

By parts (a) and (b) of \cite[Proposition 2]{RadFact} it follows that $\DS{L(x) > 1 + \frac{\alpha}{x}}$, since $x > \alpha$, and by \cite[Lemma 1]{RadFact} we have $\DS{P(2x) > 1 + \frac{\ln(2x)}{2x}}$, since $2x > 1$. One can easily show that
$$
L(x)P(2x) - 1 > \frac{1}{x}\left(\alpha + \left(x + \alpha\right)\frac{\ln(2x)}{2x}\right)
$$
which means that $a_{\mathcal L}(x) > \ln(2x)QL(x)$,
where
$$
QL(x) = \frac{\alpha}{\ln(2x)} + \frac{1}{2}\left(1 + \frac{\alpha}{x}\right).
$$
It is clear that $QL(x)$ is a positive valued strictly decreasing function, since $2x > 1$, and that $\DS{\lim_{x \rTo \infty}QL(x) = \frac{1}{2}}$. We have shown that
\begin{equation}\label{EqFFYF}
\ln(2x)QL(x) < a_{\mathcal L}(x), a_{\mathcal R}(x) < \ln(2x)QR(x)
\end{equation}

Note that (\ref{EqDD}) and (\ref{EqFFYF}) imply
\begin{equation}\label{EqFFYZ}
\ln (2n)QL(n) < \sigma_n < \ln(2n)QR(n) + 1
\end{equation}
for all $n \geq 2$. We consider the difference between these estimates for $\sigma_n$ by examining the difference between $\ln(2x)QR(x)$ and $\ln(2x)QR(x)$. Let $D(x) = \ln(2x)QR(x) - \ln(2x)QL(x)$. Since
\begin{equation}\label{EqGG}
\ln(2x)QR(x) = \frac{2x+ x\ln(2x)}{2x - \ln(2x)} = 1 + \frac{(1+x)\ln(2x)}{2x - \ln(2x)}
\end{equation}
and
\begin{equation}\label{EqHH}
\ln(2x)QL(x) = \alpha + \ln(2x)\left(\frac{x + \alpha}{2x}\right),
\end{equation}
it follows that
$$
D(x) =  (1 - \alpha)\left(\frac{2x}{2x - \ln(2x)}\right) + \left(\frac{x + \alpha}{2x - \ln(2x)}\right) \left(\frac{\ln(2x)}{\sqrt{2x}}\right)^2.
$$
Since $\DS{\lim_{x \rTo \infty}\frac{\ln x}{x^a} = 0}$ for all $a > 0$, we have $\DS{\lim_{x \rTo \infty}D(x) = 1 - \frac{\ln \pi}{2}}$. Observe that the parenthesized quotients in the preceding equation for $D(x)$ are all positive valued functions on $[e/2, \infty)$, the first two strictly decreasing there and the third on $[e^2/2, \infty)$, and finally that $1 - \alpha > 0$. Collecting results:
\begin{Prop}\label{PropDXResults}
Let $QL(x)$, $QR(x)$, and $D(x)$ be defined on $[e/2, \infty)$ as above. They are positive valued and:
\begin{enumerate}
\item[{\rm (a)}] $QL(x)$ and $QR(x)$ strictly decreasing on $[e/2, \infty)$ and $D(x)$ is strictly decreasing on $[e^2/2, \infty)$.
\item[{\rm (b)}] $\DS{\lim_{x \rTo \infty}QL(x) = \lim_{x \rTo \infty}QR(x) = \frac{1}{2}}$
\item[{\rm (c)}] $\DS{\lim_{x \rTo \infty}D(x) = 1 - \frac{\ln \pi}{2}}$.
\end{enumerate} \qed
\end{Prop}
\medskip

\begin{Cor}\label{CorollaryDXResults}
For $D(x)$, $QL(n)$, and $QR(n)$:
\begin{enumerate}
\item[{\rm (a)}] $D(x) < 1$ for all $x \geq 3.92466$.
\item[{\rm (b)}] Let $n \geq 4$. Then there are at most two integers $\ell$ in the open interval $(\ln(2n)QL(n), \; \ln(2n)QR(n)+ 1)$ and the integer $\sigma_n$ is such an $\ell$.
\end{enumerate} \qed
\end{Cor}

\pf
One can show $D(3.92465) > 1 > D(3.92466)$ and $3.92466 > e^2/2$. Thus part (a) follows by parts (a) and (c) Proposition \ref{PropDXResults}. As for part (b), note that the difference between $\ln(2n)QR(n) + 1$ and $\ln(2n)QL(n)$ is $D(n) + 1$. Thus if $n \geq 4$ then $D(n) + 1 < 2$ by part (a).  At this point (\ref{EqFFYF}) completes the proof of part (b).
\qed
\medskip

Apropos of part (a) of Corollary \ref{CorollaryDXResults}, as $0.57236 < 1 - \alpha < 0.57237$, it follows that $1.57236 < D(x) + 1 < 2$ for all $x \geq 3.92466$. Thus the length of the interval of part (b) is greater than $1.57236$ and less than $2$. By part (a) of Proposition \ref{PropDXResults}, the length of the interval $[a_{{\mathcal L}(n)}, \; a_{{\mathcal R}(n)} + 1]$ is greater than $1$. This interval, implicit in (\ref{EqDD}), is contained in the one of part (b) by (\ref{EqFFYF}); thus its length is less than $2$ when $n \geq 4$.

We examine how frequently increase in term value occurs in the sequence $\sigma_1, \sigma_2, \sigma_3, \ldots \;$. Recall its properties described in (\ref{EqB}) and (\ref{EqC}). Note that value increase occurs infinitely often since $\lim_{x \rTo \infty} \sigma_n = \infty$ by (\ref{EqFFYZ}) and that $\DS{\lim_{x \rTo \infty}QL(x) = \frac{1}{2}}$.

Let $i \geq 1$. Set $n(i) = n_{i+1} - n_i$. Then
$$
\sigma_{n_{i+1} + 1} = \sigma_{n_{i+1}+1} = \sigma_{n_{i+1}} + 1 = \sigma_{n_i} + 2.
$$
Therefore
$$
\ln(2n_i)QL(n_i)+ 2 < \sigma_{n_i} + 2 < \ln(2(n_{i+1} + 1))QR(n_{i+1} + 1) + 1
$$
by (\ref{EqFFYZ}). This means
\begin{equation}\label{EqZ}
\ln(2n_i)QL(n_i)+ 1 < \ln(2(n_{i+1} + 1))QR(n_{i+1} + 1)
\end{equation}
or
\begin{equation}\label{EqZZ}
\ln(2n_i)QL(n_i) < F(n_{i+1} + 1),
\end{equation}
where
$$
F(x) = (1+x)\left(\frac{\ln (2x)}{2x - \ln(2x)}\right)
$$
for $x > 0$, by (\ref{EqGG}).

At this point we assume $x \geq 1/2$. When $y = \ln (2x)$ then $y \geq 0$ and $e^y = 2x$. From this point on we assume $y \geq 0$. Let
$$
G(y) = \left(1 + \frac{e^y}{2}\right)\left(\frac{y}{e^y - y}\right) = \frac{1}{2}\left(\frac{2y + ye^y}{e^y - y}\right).
$$
Then $F(x) = G(\ln(2x))$. Observe that
\begin{equation}\label{EqZP}
\lim_{y \rTo \infty} (G(y) - \frac{1}{2}y) = 0
\end{equation}
and
$$
G'(y) = \frac{e^yGN(y)}{2(e^y - y)^2},
$$
where $GN(y) = e^y - (y^2 + 2y -2)$. By considering the first, second, and third derivatives of $GN(y)$ one can show that $GN(y)$ is strictly increasing on $[y_0, \infty)$, where $1.67834 < y_0 < 1.67845 = r$. As $QN(r) > 0$, it follows that $F(x)$ is strictly increasing on $[e^r/2,\infty)$. Note that $2.6784 > e^r/2$.

We return to (\ref{EqZZ}) and consider $\ln(2x)QL(x)$ in terms of $y = \ln(2x)$ as we did $F(x)$. Let
$$
H(y) = \alpha\left(1 + \frac{y}{e^y}\right) + \frac{y}{2} = \frac{1}{2}\left(\frac{(\ln\pi)(y + e^y) + ye^y}{e^y}\right).
$$
Then $\ln(2x)QL(x) = H(\ln(2x))$. Observe that
\begin{equation}\label{EqZQ}
\lim_{y \rTo \infty} (H(y) - \frac{1}{2}(y + \ln \pi)) = 0
\end{equation}
and
$$
H'(y) =  \frac{1}{2}\left(\frac{e^y - y + 1}{e^y}\right).
$$
It is easy to see that $H(y)$ is strictly increasing on $[0, \infty)$.

Observe that
\begin{equation}\label{EqZW}
\lim_{y \rTo \infty}(G(y+a) - H(y)) = \DS{\frac{1}{2}\left(a - \ln \pi\right)})
\end{equation}
follows by (\ref{EqZP}) and (\ref{EqZQ}). We examine (\ref{EqZW}) in detail when $a \geq 0$.

Let $z = y + a$ as a matter of convenience for calculations. Then $e^z = Ae^y$, where $A = e^a$. It is not hard to see that
$$
G(z) - H(y) = \frac{1}{2}\left(\frac{A(a - \ln\pi)e^{2y} + B(a, y)e^y + C(a,y)}{(e^z - z)e^y}\right),
$$
where
$$
B(a, y) = y^2 + (2 + a + (\ln\pi)(1 - A))y + a(2 + \ln\pi)
$$
and
$$
C(a, y) = y(y+a)(\ln\pi).
$$
Since $B(a, y)$ and $C(a, y)$ are quadratics whose graphs open upward, for some $c > 0$ it follows that $B(a, y)$, $C(a, y) > 0$ for all $y \geq c$. Hence if $a \geq \ln \pi$ then $G(y+a) - H(y) > 0$ for $y \geq c$.

Suppose $a < \ln \pi$. Then $G(y + a) - H(y) < 0$ eventually. Observe that $G(y + a) - H(y) < 0$ if and only if
\begin{equation}\label{EqZXZX}
\frac{B(a, y)}{e^y} + \frac{C(a,y)}{e^{2y}} < A(\ln \pi - a).
\end{equation}
Now $0 \leq a < \ln \pi$ by assumption. Therefore $1 \leq A < \pi$. We replace $B(a, y)$ and $C(a, y)$ by over estimates $B_0(y) = y^2 + (2 + \ln \pi)y + (\ln\pi)(2 + \ln\pi)$ and $C_0(y) = y(y + \ln\pi)(\ln\pi)$ respectively. Therefore(\ref{EqZXZX}) holds for $y$ if
\begin{equation}\label{EqZXZXX}
\frac{B_0(y)}{e^y} + \frac{C_0(y)}{e^{2y}} < A(\ln \pi - a).
\end{equation}
Let $\ell$ be a non-negative and let $m$ be a positive integer. Then $\DS{\frac{y^\ell}{e^{my}}}$ is strictly decreasing on $\DS{[\frac{\ell}{m}, \infty)}$. Hence the left hand side of the inequality of (\ref{EqZXZXX}) is a strictly decreasing function on $[2, \infty)$. Therefore if $c \geq 2$ and (\ref{EqZXZXX}) holds for $y = c$ then (\ref{EqZXZX}) holds for all $y \geq c$.

Let $a = \ln 3$. Then $0 < a < \ln\pi$. One can show that (\ref{EqZXZXX}) does not hold for $y = 6.06520$ and does hold for $y = 6.06521$. Consequently (\ref{EqZXZX}) holds for all $y \geq 6.06521$ when $a = \ln 3$; hence
\begin{equation}\label{EqZXZXXX}
G(y + \ln 3) < H(y)
\end{equation}
for all $y \geq 6.06521$ and, as $215.30654 < e^{6.06521}/2 < 215.30655$, we have shown that
\begin{equation}\label{EqZXZXXXZ}
F(3x) < \ln(2x)QL(x)
\end{equation}
for all $x \geq 215.30655$.
Therefore
\begin{equation}\label{EqZXZXXXWS}
F(3n_i) < \ln(2n_i)QL(n_i) < F(n_{i+1} + 1)
\end{equation}
for all $n_i \geq 216$. Assume $n_i$ satisfies the preceding inequality. Since $F(x)$ is strictly increasing on $[2.7684, \infty)$, it follows that $3n_i < n_{i+1} + 1$ and therefore $3n_i \leq n_{i+1}$.

Using (\ref{EqA}) one can show that $\sigma_1 = \cdots = \sigma_3 = 2$, $\sigma_4 = \cdots = \sigma_{54} = 3$, $\sigma_{55} = \cdots = \sigma_{458} = 4$, $\sigma_{459} = \cdots = \sigma_{3480} = 5$, and $\sigma_{3481} = \cdots = \sigma_{25867} = 6$, and $\sigma_{25868} = \cdots = \sigma_{191351} = 7$, and $\sigma_{191352} = 8$. Thus $n_1 = 3$, $n_2 = 54$, $n_3 = 458$, $n_4 = 3480$, $n_5 = 25867$, and $n_6 = 191351$. Observe that

We have shown in particular that:
\begin{Cor}\label{CorollaryNI}
$3n_i \leq n_{i+1}$, in particular $\DS{\frac{n_{i+1}}{n_i} \geq 3}$, for all $i \geq 1$. \qed
\end{Cor}
\medskip

There is an interesting relationship between \cite[Proposition 3]{RadFact} and the corollary. Let $n \geq 1$ and $S_n = T_{n+1} - T_n$, where $T_n = e\sqrt[n]{n!}\;$ is defined just after (\ref{EqA}). Write $\DS{S_n = 1 + \frac{a(n)}{n}}$. Then $\DS{\lim_{n \rTo \infty} a(n) = \frac{1}{2}}$ by \cite[Theorem 1(d)]{RadGamma}. Let $\DS{a = \frac{1}{2}}$. Suppose that $\DS{S_n < 1 + \frac{a}{n}}$ for all $n \geq 1$. If this is the case then \cite[Proposition 3]{RadFact} implies Corollary \ref{CorollaryNI} as well.

We end this paper with comments about the sequence $n_1, n_2, n_3, \ldots\;\;$. First of all
\begin{equation}\label{EqSIG}
\sigma_{n_\ell}= \ell + 1
\end{equation}
for all $\ell \geq 1$, which follows by induction on $\ell$. Suppose that $n \geq 4$ and set $\sigma_n = \ell + 1$.  Then $\sigma_n = \sigma_{n_\ell}$ by (\ref{EqSIG}). Thus $\ell \geq 2$, as $n \geq 4$, and therefore $n > n_{\ell - 1}$ as $\sigma_n > \sigma_{n_{\ell - 1}}$. On closer scrutiny the quotients $\DS{\frac{n_{i+1}}{n_i}}$ seem to form a nice pattern. A few calculations
$$
\DS{\frac{n_2}{n_1} = 18}, \;\; \DS{8.48 < \frac{n_3}{n_2} < 8.49}, \;\; \DS{7.59 < \frac{n_4}{n_3} < 7.6},
$$
$$
\DS{7.43 < \frac{n_5}{n_4} < 7.44},\;\; \mbox{and} \;\;\DS{7.39 < \frac{n_6}{n_5} < 7.40}
$$
suggest
\begin{Conj}
The sequence $\DS{\frac{n_2}{n_1}, \frac{n_3}{n_2}, \frac{n_4}{n_3}, \ldots \;}$ of quotients is strictly decreasing and $\DS{\lim_{n \rTo \infty} \frac{n_{i+1}}{n_i} = e^2}$.
\end{Conj}
%Indeed computation of the quotients through $i = 50$ bears this out.

\end{document}